\documentclass{amsart}
\usepackage{mathrsfs}
\usepackage{amsfonts}
\usepackage{amssymb}
\usepackage[papersize={7.3in,10in},textwidth=13.4cm,textheight=20.6cm,centering]{geometry}

\usepackage[dvipdfm,colorlinks=true,citecolor=blue,linkcolor=blue]{hyperref}
\hypersetup{
pdfstartpage=1,
pdfstartview=FitH,}

\allowdisplaybreaks

\newtheorem*{theorem*}{Theorem}
\newtheorem*{corollary*}{Corollary}
\newtheorem{theorem}{Theorem}

\theoremstyle{remark}

\newcommand{\ud}{\mathrm{d}}

\newcommand{\e}{\mathbf{e}}

\title[Quadratic residues and non-residues for Piatetski-Shapiro primes]{Quadratic residues and non-residues for infinitely many Piatetski-Shapiro primes}

\author{Ping Xi}

\address{School of Mathematics and Statistics, Xi'an Jiaotong University, Xi'an 710049, P.R. China}

\email{xprime@163.com}



\begin{document}

\begin{abstract}In this paper, we prove a quantitative version of the statement that every nonempty finite subset of $\mathbb{N}^+$ is a set of quadratic residues for infinitely many primes of the form $[n^c]$ with $1\leqslant c\leqslant243/205$. Correspondingly, we can obtain a similar result for the case of quadratic non-residues under reasonable assumptions. These results generalize the previous ones obtained by S. Wright in certain aspects.\end{abstract}
\maketitle

\section{Introduction}
Distributions of quadratic residues and non-residues have received great attention over the recent decades. S. Wright \cite{W1} has proved by combinatorial and reasonably elementary method that

(1) {\it Every nonempty finite subset of $\mathbb{N}^+$ is a set of quadratic residues for infinitely
many primes. }

(2) {\it For every nonempty finite subset $\mathcal{A}\subseteq\mathbb{N}^+$, which contains no squares, we write $\Pi_{\text{odd}}(a)=\{p \text{ prime}:p^\alpha\|a,\alpha\text{ {\rm odd}}\}$, $\Pi=\bigcup_{a\in\mathcal{A}}\Pi_{\text{{\rm odd}}}(a)$ and $\mathfrak{P}=\{\Pi_{\text{{\rm odd}}}(a):a\in\mathcal{A}\}$.

If $\mathfrak{P}$ contains at most $4$ elements, then $\mathcal{A}$ is a set of quadratic non-residues for infinitely
many primes if and only if $\mathfrak{P}$ does not contain a $3$-cycle included in $\{\mathcal{C}:\mathcal{C}\neq\emptyset,
\mathcal{C}\subseteq\Pi\}$.}

We intend to give a quantitative version of such results in the present paper. To be precise, for an arbitrarily given set $\mathcal{S}$ of finite cardinality $S=\#\mathcal{S}$, consider the following quantity
\[\mathcal{S}(x)=\#\left\{p\in\mathcal{P}\cap(x,2x]:\left(\frac{s}{p}\right)=1\text{ for each }s\in\mathcal{S}\right\},\]
and correspondingly
\[\overline{\mathcal{S}}(x)=\#\left\{p\in\mathcal{P}\cap(x,2x]:\left(\frac{s}{p}\right)=-1\text{ for each }s\in\mathcal{S}\right\},\]
where $\mathcal{P}$ be the set consisting of all the primes, $(\frac{\cdot}{p})$ is the Legendre symbol mod $p$. Clearly, we aim to show that $\mathcal{S}(x)\gg x^{\theta_1}$ and $\overline{\mathcal{S}}(x)\gg x^{\theta_2}$ for certain $\theta_1,\theta_2>0$ under reasonable assumptions.

On the other hand, a classical result due to I.I. Piatetski-Shapiro states that there exists $c>1$ such that there are infinitely many
primes of the form $[n^c]$ with $n\in\mathbb{N}^+$, where $[x]$ denotes the integral part of $x$. More precisely,
Piatetski-Shapiro \cite{PS} proved that
\[\sum_{\substack{n\leqslant x\\ [n^c]\in\mathcal{P}}}1=\left(\frac{1}{c}+o(1)\right)\frac{x}{\log x}, \text{  as }x\rightarrow+\infty\]
for $1<c<12/11=1.\overline{09}$. There are many works devoted to enlarge the reasonable range of $c$, see \cite{HB,J1,J2,Ko,Ku,LR,RW} for instance. The most recent result is due to J. Rivat and J. Wu \cite{RW}, who succeeded in showing that
\[\sum_{\substack{n\leqslant x\\ [n^c]\in\mathcal{P}}}1\gg\frac{x}{c\log x}, \text{  as }x\rightarrow+\infty\]
for the wider range $1\leqslant c\leqslant243/205=1.18536\ldots.$ Moreover, there are many other problems concerning primes, as well as some related to the analytic theory of automorphic forms (see \cite{BZ} for instance), which also restrict the primes to be of the form $[n^c]$.

Motivated by these works on the distribution of Piatetski-Shapiro primes, we can consider
the following quantity
\[\mathcal{S}_c(x)=\#\left\{x<n\leqslant2x:[n^c]\in\mathcal{P},\text{ and}\left(\frac{s}{[n^c]}\right)=1\text{ for each }s\in\mathcal{S}\right\},\]and correspondingly
\[\overline{\mathcal{S}}_c(x)=\#\left\{x<n\leqslant2x:[n^c]\in\mathcal{P},\text{ and}\left(\frac{s}{[n^c]}\right)=-1\text{ for each }s\in\mathcal{S}\right\}.\]
Clearly, $\mathcal{S}_1(x)=\mathcal{S}(x)$, $\overline{\mathcal{S}}_1(x)=\overline{\mathcal{S}}(x)$.

We shall prove by analytic methods that
\begin{theorem}\label{thm:1} Every nonempty finite set $\mathcal{S}\subseteq\mathbb{N}^+$ is a set of quadratic residues of infinitely many primes of the form
$[n^c]$ with $1\leqslant c\leqslant243/205.$

To be precise, let $c$ be a number with $1\leqslant c\leqslant243/205.$ Then for sufficiently large $x$, we have
\begin{align*}\mathcal{S}_c(x)=\bigg(\frac{\#\mathfrak{H}+1}{c\cdot2^{\#\mathcal{S}}}+o(1)\bigg)\frac{x}{\log x},\end{align*}
where $\mathfrak{H}$ denotes the collection of all the nonempty subsets of $\mathcal{S}$ satisfying that the product of all the elements of each subset is a square integer.
\end{theorem}

\begin{theorem}\label{thm:2}Consider all of the nonempty subsets of the nonempty finite set $\mathcal{S}\subseteq\mathbb{N}^+,$ the product of the elements of each one is a square. If such subsets of even cardinalities
are not less than those of odd cardinalities, then
$\mathcal{S}$ is a set of quadratic non-residues for infinitely
many primes of the form $[n^c]$ with $1\leqslant c\leqslant243/205.$

More precisely, we have
\begin{align*}\overline{\mathcal{S}}_c(x)&=\frac{1}{c\cdot2^{\#\mathcal{S}}}\bigg(1+\sum_{\mathcal{T}\in\mathfrak{H}}(-1)^{\#\mathcal{T}}+o(1)\bigg)\frac{x}{\log x},\end{align*} where $\mathfrak{H}$ is the same with that in Theorem $\ref{thm:1}$.
\end{theorem}

\begin{corollary*}For any given finitely many primes, all of them are the quadratic residues for certain infinitely many primes of the form $[n^c]$ with $1\leqslant c\leqslant243/205,$ and in the meanwhile, are also the quadratic non-residues for infinitely many primes of the form $[n^c]$ with $1\leqslant c\leqslant243/205.$
\end{corollary*}

In fact, during the following procedure of the proof, we have applied a well-known result
\begin{align*}\lim_{x\rightarrow+\infty}\frac{1}{\pi(x)}\sum_{p\leqslant x}\bigg(\frac{a}{p}\bigg)=0\end{align*}
for any $a$ which is not a square. Let $\mathcal{P}_0\subseteq\mathcal{P}$ be an infinite set. If the distribution of the primes in $\mathcal{P}_0$ is also well enough, such that the characteristic function of $\mathcal{P}_0$ is asymptotically orthogonal to the quadratic characters, i.e.,
\begin{align*}\lim_{x\rightarrow+\infty}\frac{1}{\#\mathcal{P}_0\cap[1,x]}\sum_{p\in\mathcal{P}_0\cap[1,x]}\bigg(\frac{a}{p}\bigg)=0\end{align*}for any non-square integer $a$, then the results in Theorems \ref{thm:1} and \ref{thm:2} will also hold if we restrict those primes to be in $\mathcal{P}_0$ instead of to be only of the form $[n^c].$

\bigskip

\noindent \textbf{Notation.} As usual, $\mathcal{P}$ denotes the set consisting of all the primes and $p\in\mathcal{P}$. $\#\mathcal{A}$ stands for the cardinality of the set $\mathcal{A}$. $[x]$ denotes the largest integer not exceeding $x$ and $\psi(x)=x-[x]-1/2.$ $\Lambda(n)$ and $\mu(n)$ denote the von Mangoldt function and M\"{o}bius function, respectively. $(\frac{\cdot}{n})$ denotes the Legendre-Jacobi symbol of an odd integer $n$. We denote by $\square$ the square of an integer. For any set $\mathcal{B}$, we write $(p,\mathcal{B})=1$, yielding that $(p,b)=1$ for each $b\in \mathcal{B}$. For convenience, we write $\gamma=1/c$ throughout this paper.

\bigskip

\section{Proof of Theorem \ref{thm:1}}

First, by virtue of Legendre symbol, we can write
\[\mathcal{S}_c(x)=\frac{1}{2^S}\sum_{\substack{x<n\leqslant2x\\ [n^c]\in\mathcal{P}\\([n^c],\mathcal{S})=1}}\prod_{s\in\mathcal{S}}\bigg(1+\bigg(\frac{s}{[n^c]}\bigg)\bigg).\]In fact, we can take $x$ sufficiently large, such that the condition $([n^c],\mathcal{S})=1$ can be omitted, for instance we can take $x>\prod_{s\in\mathcal{S}}s^\gamma$.

For any prime $p$, the condition $[n^c]=p$ is equivalent to $p\leqslant n^c<p+1,$ i.e., \[-(p+1)^\gamma<-n\leqslant-p^\gamma,\] we have
\[\mathcal{S}_c(x)=\frac{1}{2^S}\sum_{x^c<p\leqslant (2x)^c}([-p^\gamma]-[-(p+1)^\gamma])\prod_{s\in\mathcal{S}}\bigg(1+\bigg(\frac{s}{p}\bigg)\bigg).\]

Observing that $[x_1]-[x_2]=(x_1-x_2)-(\psi(x_1)-\psi(x_2))$, thus
we only require to study the quantities
\[\Sigma_1(N)=\sum_{N<p\leqslant\delta N}((p+1)^\gamma-p^\gamma)\prod_{s\in\mathcal{S}}\bigg(1+\bigg(\frac{s}{p}\bigg)\bigg),\]
and
\[\Sigma_2(N)=\sum_{N<p\leqslant\delta N}(\psi(-p^\gamma)-\psi(-(p+1)^\gamma))\prod_{s\in\mathcal{S}}\bigg(1+\bigg(\frac{s}{p}\bigg)\bigg).\]

\subsection{\bf Estimate for $\Sigma_1(N)$}

First, we have
\begin{align*}\prod_{s\in\mathcal{S}}\bigg(1+\bigg(\frac{s}{p}\bigg)\bigg)=1+\sum_{\emptyset\neq\mathcal{T}\subseteq\mathcal{S}}\prod_{s\in\mathcal{T}}\bigg(\frac{s}{p}\bigg).\end{align*}

For certain $\emptyset\neq\mathcal{T}\subseteq\mathcal{S}$, the product of all the elements of $\mathcal{T}$ is the square of certain integer, namely
\[\prod_{s\in\mathcal{T}}s=\square.\]We can collect these subsets $\mathcal{T}\subseteq\mathcal{S}$ of such a property, to construct a set $\mathfrak{H}=\mathfrak{H}(\mathcal{S})$, saying
\[\mathfrak{H}=\left\{\mathcal{T}\subseteq\mathcal{S}:\mathcal{T}\neq\emptyset,\prod_{s\in\mathcal{T}}s=\square\right\}.\]
Now we have
\begin{align}\label{eq:1}\prod_{s\in\mathcal{S}}\bigg(1+\bigg(\frac{s}{p}\bigg)\bigg)
=\#\mathfrak{H}+1+\sum_{\substack{\emptyset\neq\mathcal{T}\subseteq\mathcal{S}\\ \mathcal{T}\not\in\mathfrak{H}}}\prod_{s\in\mathcal{T}}\bigg(\frac{s}{p}\bigg).\end{align}

The first two terms in (\ref{eq:1}) contribute to $\Sigma_1(N)$ as
\begin{align*}\Sigma_{11}(N)&=(\#\mathfrak{H}+1)\sum_{N<p\leqslant\delta N}((p+1)^\gamma-p^\gamma)\\
&=(\#\mathfrak{H}+1)\int_N^{\delta N}((y+1)^\gamma-y^\gamma)\ud\pi(y)\\
&=(\#\mathfrak{H}+1+o(1))\frac{N^\gamma}{\log N}.\end{align*}

Correspondingly, the third term in (\ref{eq:1}) gives a contribution to $\Sigma_1(N)$ as
\[\Sigma_{12}(N)=\sum_{N<p\leqslant\delta N}((p+1)^\gamma-p^\gamma)\sum_{\substack{\emptyset\neq\mathcal{T}\subseteq\mathcal{S}\\ \mathcal{T}\not\in\mathfrak{H}}}\prod_{s\in\mathcal{T}}\bigg(\frac{s}{p}\bigg).\]

It is well-known that
\begin{align*}\sum_{p\leqslant N}\bigg(\frac{a}{p}\bigg)=o(\pi(N)),\ \text{for any }a\neq\square,\end{align*}
thus by patrial summation we get
\begin{align*}\Sigma_{12}(N)&=\int_N^{\delta N}((y+1)^\gamma-y^\gamma)\ud\bigg(\sum_{N<p\leqslant y}\sum_{\substack{\emptyset\neq\mathcal{T}\subseteq\mathcal{S}\\ \mathcal{T}\not\in\mathfrak{H}}}\prod_{s\in\mathcal{T}}\bigg(\frac{s}{p}\bigg)\bigg)
=o\left(\frac{N^\gamma}{\log N}\right).\end{align*}

Hence we find that
\begin{align}\label{eq:2}\Sigma_1(N)=(\delta^\gamma-1)(\#\mathfrak{H}+1+o(1))\frac{N^\gamma}{\log N}.\end{align}

\subsection{\bf Estimate for $\Sigma_2(N)$}

First, we have
\begin{align*}\Sigma_2(N)=\sum_{N<p\leqslant\delta N}(\psi(-p^\gamma)-\psi(-(p+1)^\gamma))\bigg(\#\mathfrak{H}+1+\sum_{\substack{\emptyset\neq\mathcal{T}\subseteq\mathcal{S}\\ \mathcal{T}\not\in\mathfrak{H}}}\prod_{s\in\mathcal{T}}\bigg(\frac{s}{p}\bigg)\bigg).\end{align*}

Note that we have the truncated Fourier expansion (see \cite{GK}, Section 4.6) for $\psi(x)$
\begin{align}\label{eq:3}\psi(x)=\psi^*(x)+O(\delta(x)),\end{align}where
\[\psi^*(x)=\sum_{0<|j|\leqslant J}a(j)e(jx),\ \ \ \delta(x)=\sum_{0\leqslant|j|\leqslant J}b(j)e(jx),\]and $a(j)\ll|j|^{-1},b(j)\ll J^{-1}$. Moreover, $\delta(x)$ is a non-negative function. Here $J$ is a parameter to be chosen later according to the particular properties of the objects in question. Now we have
\begin{align*}\Sigma_2(N)&=\sum_{N<p\leqslant\delta N}(\psi^*(-p^\gamma)-\psi^*(-(p+1)^\gamma))\bigg(\#\mathfrak{H}+1+\sum_{\substack{\emptyset\neq\mathcal{T}\subseteq\mathcal{S}\\ \mathcal{T}\not\in\mathfrak{H}}}\prod_{s\in\mathcal{T}}\bigg(\frac{s}{p}\bigg)\bigg)\\
&\ \ \ \ \ \ +O\bigg(\sum_{N<p\leqslant\delta N}(\delta(-p^\gamma)+\delta(-(p+1)^\gamma))\bigg).\end{align*}

It has been shown in \cite{GK} that
\begin{align*}\sum_{N<n\leqslant\delta N}\delta(-n^\gamma)\ll N^\gamma\log^{-2}N+N^{1/2}\log N,\end{align*}
thus for $\gamma>1/2,$ we have
\[\Sigma_2(N)=(\#\mathfrak{H}+1)\Sigma_{21}(N)+\Sigma_{22}(N)+o\left(\frac{N^\gamma}{\log N}\right),\]where
\begin{align*}\Sigma_{21}(N)&=\sum_{N<p\leqslant\delta N}(\psi^*(-p^\gamma)-\psi^*(-(p+1)^\gamma)),\end{align*}
and
\begin{align*}\Sigma_{22}(N)&=\sum_{N<p\leqslant\delta N}(\psi^*(-p^\gamma)-\psi^*(-(p+1)^\gamma))\sum_{\substack{\emptyset\neq\mathcal{T}\subseteq\mathcal{S}\\ \mathcal{T}\not\in\mathfrak{H}}}\prod_{s\in\mathcal{T}}\bigg(\frac{s}{p}\bigg).\end{align*}

Next, we aim to show that
\begin{align*}|\Sigma_{21}(N)|+|\Sigma_{22}(N)|=o\left(\frac{N^\gamma}{\log N}\right).\end{align*} To this end, we just need to prove the following estimates
\begin{align*}L_1(N)&:=\sum_{N<n\leqslant M}\Lambda(n)(\psi^*(-n^\gamma)-\psi^*(-(n+1)^\gamma))=o(N^\gamma),\end{align*}
and
\begin{align*}L_2(N)&:=\sum_{N<n\leqslant M}\Lambda(n)(\psi^*(-n^\gamma)-\psi^*(-(n+1)^\gamma))\bigg(\frac{s}{n}\bigg)=o(N^\gamma)\end{align*}
for any $N<M\leqslant\delta N$ and squarefree $s$. Here $(\frac{\cdot}{n})$ should be regarded as the Jacobi symbol and the summations are both restricted to be over odd integers $n$.

For $L_1(N)$, such an estimate has been achieved for $205/243\leqslant\gamma\leqslant1$ in \cite{RW}. Hence we only need to consider $L_2(N)$.

In view of (\ref{eq:3}), and writing
\[\phi_j(x)=1-e(j((x+1)^\gamma-x^\gamma)),\]
we have
\begin{align*}L_2(N)&=\sum_{0<|j|\leqslant J}a(j)\sum_{N<n\leqslant M}\Lambda(n)(e(jn^\gamma)-e(j(n+1)^\gamma))\bigg(\frac{s}{n}\bigg)\\
&=\sum_{0<|j|\leqslant J}a(j)\sum_{N<n\leqslant M}\Lambda(n)e(jn^\gamma)\phi_j(n)\bigg(\frac{s}{n}\bigg)\\
&\ll\sum_{j\leqslant J}j^{-1}\left|\sum_{N<n\leqslant M}\Lambda(n)e(jn^\gamma)\phi_j(n)\bigg(\frac{s}{n}\bigg)\right|.\end{align*}

By partial summation, we get
\begin{align*}L_2(N)&\ll\sum_{j\leqslant J}j^{-1}\left|\phi_j(M)\sum_{N<n\leqslant M}\Lambda(n)e(jn^\gamma)\bigg(\frac{s}{n}\bigg)\right|\\
&\ \ \ \ \ \ +\int_N^{M}\sum_{j\leqslant J}j^{-1}\left|\frac{\partial\phi_j(t)}{\partial t}\sum_{N<n\leqslant t}\Lambda(n)e(jn^\gamma)\bigg(\frac{s}{n}\bigg)\right|\ud t.\end{align*}
Observing that
\[\phi_j(t)\ll jN^{\gamma-1},\ \ \frac{\partial\phi_j(t)}{\partial t}\ll jN^{\gamma-2},\]
thus there must exist certain $N'$ with $N<N'\leqslant\delta N$, such that
\begin{align*}L_2(N)\ll N^{\gamma-1}\sum_{j\leqslant J}\left|\sum_{N<n\leqslant N'}\Lambda(n)e(jn^\gamma)\bigg(\frac{s}{n}\bigg)\right|.\end{align*}

Now we can rewrite $\Lambda(n)$ in the following manner
\[\Lambda(n)=-\sum_{\substack{kl=n\\k>v\\l>u}}\Lambda(k)\sum_{\substack{d|l\\d\leqslant u}}\mu(d)+\sum_{\substack{kl=n\\l\leqslant u}}\mu(k)\log l-\sum_{\substack{klm=n\\k\leqslant v\\m\leqslant u}}\Lambda(k)\mu(m),\]where $u,v$ are two positive real numbers and $v<n$. The idea of such an identity is originally due to R.C. Vaughan \cite{V}.

With the help of Vaughan's identity, we have
\begin{align*}\sum_{N<n\leqslant N'}\Lambda(n)e(jn^\gamma)\bigg(\frac{s}{n}\bigg)
&=-\sum_{\substack{N<mn\leqslant N'\\m>u\\n>v}}\Lambda(n)a(m)\bigg(\frac{s}{mn}\bigg)e(jm^\gamma n^\gamma)\\
&\ \ \ \ \ +\sum_{\substack{N<mn\leqslant N'\\n\leqslant u}}\mu(m)\log n\bigg(\frac{s}{mn}\bigg)e(jm^\gamma n^\gamma)\\
&\ \ \ \ \ -\sum_{\substack{N<mn\leqslant N'\\m\leqslant uv}}b(m)\bigg(\frac{s}{mn}\bigg)e(jm^\gamma n^\gamma),\end{align*}
where
\[a(m)=\sum_{\substack{d|m\\d\leqslant u}}\mu(d),\ \ \ b(m)=\sum_{\substack{de=m\\d\leqslant v\\e\leqslant u}}\Lambda(d)\mu(e)\]

There is an extra factor $(\frac{s}{mn})$ compared to that in \cite{RW}. Fortunately, for any integers $m$ and $n$, we have
\[\bigg(\frac{s}{mn}\bigg)=\bigg(\frac{s}{m}\bigg)\bigg(\frac{s}{n}\bigg),\]
thus we have
\begin{align*}\sum_{N<n\leqslant N'}\Lambda(n)e(jn^\gamma)\bigg(\frac{s}{n}\bigg)
&=-\sum_{\substack{N<mn\leqslant N'\\m>u\\n>v}}a(m)\bigg(\frac{s}{m}\bigg)\Lambda(n)\bigg(\frac{s}{n}\bigg)e(jm^\gamma n^\gamma)\\
&\ \ \ \ \ +\sum_{\substack{N<mn\leqslant N'\\n\leqslant u}}\mu(m)\bigg(\frac{s}{m}\bigg)\bigg(\frac{s}{n}\bigg)(\log n)e(jm^\gamma n^\gamma)\\
&\ \ \ \ \ -\sum_{\substack{N<mn\leqslant N'\\m\leqslant uv}}b(m)\bigg(\frac{s}{m}\bigg)\bigg(\frac{s}{n}\bigg)e(jm^\gamma n^\gamma).\end{align*}
Now, following the similar arguments in Section 4.6 of \cite{GK} and those in \cite{RW}, we can deduce that
\begin{align*}L_2(N)\ll N^{\gamma-1}\sum_{j\leqslant J}\left|\sum_{N<n\leqslant N'}\Lambda(n)e(jn^\gamma)\bigg(\frac{s}{n}\bigg)\right|=o(N^\gamma)\end{align*}
for $205/243\leqslant\gamma\leqslant1$ by choosing suitable $J$.

Hence
\[\Sigma_2(N)=o\left(\frac{N^\gamma}{\log N}\right),\]which, together with (\ref{eq:2}), yields
\begin{align*}\mathcal{S}_c(x)=\frac{\#\mathfrak{H}+1}{2^Sc}(1+o(1))\frac{x}{\log x}\end{align*} for $1\leqslant c\leqslant243/205$ by taking $\delta=2^c,N=x^c.$ This completes the proof of Theorem \ref{thm:1}.

\bigskip

\section{Proof of Theorem \ref{thm:2}}

For sufficiently large $x>\prod_{s\in\mathcal{S}}s^\gamma$, we have
\begin{align*}\overline{\mathcal{S}}_c(x)&=\frac{1}{2^S}\sum_{\substack{x<n\leqslant2x\\ [n^c]\in\mathcal{P}}}\prod_{s\in\mathcal{S}}\bigg(1-\bigg(\frac{s}{[n^c]}\bigg)\bigg)\\
&=\frac{1}{2^S}\sum_{x^c<p\leqslant (2x)^c}([-p^\gamma]-[-(p+1)^\gamma])\prod_{s\in\mathcal{S}}\bigg(1-\bigg(\frac{s}{p}\bigg)\bigg).\end{align*}

Observing that
\begin{align*}\prod_{s\in\mathcal{S}}\bigg(1-\bigg(\frac{s}{p}\bigg)\bigg)=1+\sum_{\emptyset\neq\mathcal{T}\subseteq\mathcal{S}}
(-1)^{\#\mathcal{T}}\prod_{s\in\mathcal{T}}\bigg(\frac{s}{p}\bigg),\end{align*}
and following the similar arguments, we can show that
\begin{align*}\overline{\mathcal{S}}_c(x)&=\frac{1}{2^Sc}\bigg(1+\sum_{\mathcal{T}\in\mathfrak{H}}(-1)^{\#\mathcal{T}}+o(1)\bigg)\frac{x}{\log x}\end{align*}
for $1\leqslant c\leqslant243/205.$

Of course, we can find
\begin{align*}\sum_{\mathcal{T}\in\mathfrak{H}}(-1)^{\#\mathcal{T}}\geqslant-1\end{align*}holds for any given finite set $\mathcal{S}$, or else $\overline{\mathcal{S}}_c(x)\rightarrow-\infty$ as $x\rightarrow+\infty$, which is obviously absurd. However, if we expect $\mathcal{S}$ is a set of quadratic non-residues of certain infinitely many primes, we can assume
\begin{align}\label{eq:4}\sum_{\mathcal{T}\in\mathfrak{H}}(-1)^{\#\mathcal{T}}\geqslant0,\end{align}
which is equivalent to
\begin{align*}\#\{\mathcal{T}\in\mathfrak{H}:\#\mathcal{T}\text{ even}\}\geqslant\#\{\mathcal{T}\in\mathfrak{H}:\#\mathcal{T}\text{ odd}\}.\end{align*}
Now Theorem \ref{thm:2} follows immediately.

\bigskip

\section{Remarks}

For Theorem \ref{thm:2}, the only case we have left is
\begin{align*}\#\{\mathcal{T}\in\mathfrak{H}:\#\mathcal{T}\text{ odd}\}-\#\{\mathcal{T}\in\mathfrak{H}:\#\mathcal{T}\text{ even}\}=1,\end{align*}
which is equivalent to
\begin{align*}\sum_{\mathcal{T}\in\mathfrak{H}}(-1)^{\#\mathcal{T}}=-1.\end{align*}Clearly, this implies that
\begin{align*}\#\mathfrak{H}-1=\sum_{\mathcal{T}\in\mathfrak{H}}(1+(-1)^{\#\mathcal{T}})=2\#\{\mathcal{T}\in\mathfrak{H}:\#\mathcal{T}\text{ even}\},\end{align*}
which shows that $\#\mathfrak{H}$ must be odd. Hence, if $\#\mathfrak{H}$ is even, the condition (\ref{eq:4}) can never hold. However, for the case that $\#\mathfrak{H}$ is odd, we still don't know whether $\mathcal{S}$ is a set of quadratic non-residues for infinitely many primes.

We would like to mention at last that S. Wright \cite{W2} proved by combinatorial arguments a criterion both necessary and sufficient for $\mathcal{S}$ to
satisfy the following condition: if $n$ is a fixed nonnegative integer, then there exists infinitely many primes $p$
such that $\mathcal{S}$ contains exactly $n$ quadratic residues of $p$. And the case for quadratic non-residues is also investigated in \cite{W3}. It's quite interesting to prove some further related results by analytic methods.

\bigskip

\noindent\textbf{Acknowledgements.} The author would like to express his gratitude to Prof. Yuan Yi for the constant help and great encouragement, and also to the referee for pointing out the references \cite{W2} and \cite{W3}.

\bigskip

\bibliographystyle{plainnat}

\bigskip

\end{document}